\documentclass[12pt]{amsart}
\usepackage{amssymb,latexsym,comment,url}
\usepackage{graphicx}

\def\Q{{\mathbb Q}}
\newcommand\Sage{{\sc SageMath}}
\newcommand\PARI{{\tt Pari/gp}}
\newcommand\Magma{{\sc Magma}}
\newcommand\mongodb{{\tt MongoDB}}
\newcommand\pymongo{{\tt PyMongo}}
\newcommand\Python{{\tt Python}}
\newcommand\GitHub{{\tt GitHub}}
\newcommand\git{{\tt Git}}
\newcommand\flask{{\tt Flask}}
\newcommand\jinja{{\tt Jinja}}

                      {\hspace*{\fill}\nobreak$\Box$\par\medskip}
                       {\hspace*{\fill}\nobreak$\Box$\par\medskip}

\theoremstyle{definition}

\begin{document}

\title[The LMFDB project]{The L-functions and modular forms database project}
\author{John Cremona}
\address{Warwick Mathematics Institute, Zeeman Building, University of Warwick, Coventry, CV4 7AL, UK}
\date{18 September 2015, revised January 2016}

\begin{abstract}
The Langlands Programme, formulated by Robert Langlands in the 1960s
and since much developed and refined, is a web of interrelated theory
and conjectures concerning many objects in number theory, their
interconnections, and connections to other fields.  At the heart of
the Langlands Programme is the concept of an L-function.

The most famous L-function is the Riemann zeta-function, and as well
as being ubiquitous in number theory itself, L-functions have
applications in mathematical physics and cryptography. Two of the
seven Clay Mathematics Million Dollar Millennium Problems, the Riemann
Hypothesis and the Birch and Swinnerton-Dyer Conjecture, deal with
their properties. Many different mathematical objects are connected in
various ways to L-functions, but the study of those objects is highly
specialized, and most mathematicians have only a vague idea of the
objects outside their specialty and how everything is related. Helping
mathematicians to understand these connections was the motivation for
the L-functions and Modular Forms Database (LMFDB) project. Its
mission is to chart the landscape of L-functions and modular forms in
a systematic, comprehensive and concrete fashion. This involves
developing their theory, creating and improving algorithms for
computing and classifying them, and hence discovering new properties
of these functions, and testing fundamental conjectures.

In the lecture I gave a very brief introduction to L-functions for
non-experts, and explained and demonstrated how the large collection
of data in the LMFDB is organized and displayed, showing the
interrelations between linked objects, through our website
\url{www.lmfdb.org}. I also showed how this has been created by a
world-wide open source collaboration, which we hope may become a model
for others.
\end{abstract}

\maketitle

\section{What is the LMFDB?}

Since the early days of using computers in number theory, computations
and tables have played an important part in experimentation, for the
purpose of formulating and proving (or disproving) conjectures.

Until the World Wide Web, such tables were hard to use, let alone to
make, as they were only available in printed form, or on microfiche!
An example relevant for the LMFDB is the 1976 Antwerp IV tables of
elliptic curves, published as part of a conference proceedings in
Springer Lecture Notes in Mathematics 476, as a computer printout with
manual amendments and diagrams.

However, even since the WWW, tables and databases have been scattered
among a variety of personal web pages (including my own
\cite{ecdata}).  To use them, you had to know who to ask, download
data, and deal with a wide variety of formats.  A few had more
sophisticated interfaces, but there was no consistency.

In some areas of number theory, such as elliptic curves, the situation
is now much better and easier: packages such as \Sage\ \cite{Sage},
\Magma\ \cite{Magma} and \PARI\ \cite{PARI2} contain elliptic curve
databases (sometimes as optional add-ons, as they are large).  Also,
the internet makes accessing even ``printed'' tables much easier.  But
the data are still very scattered and incomplete.

The situation is now very much better: we have the \textbf{LMFDB}!
\medskip

\hbox to \hsize{
\hfill \includegraphics[height=3.2in]{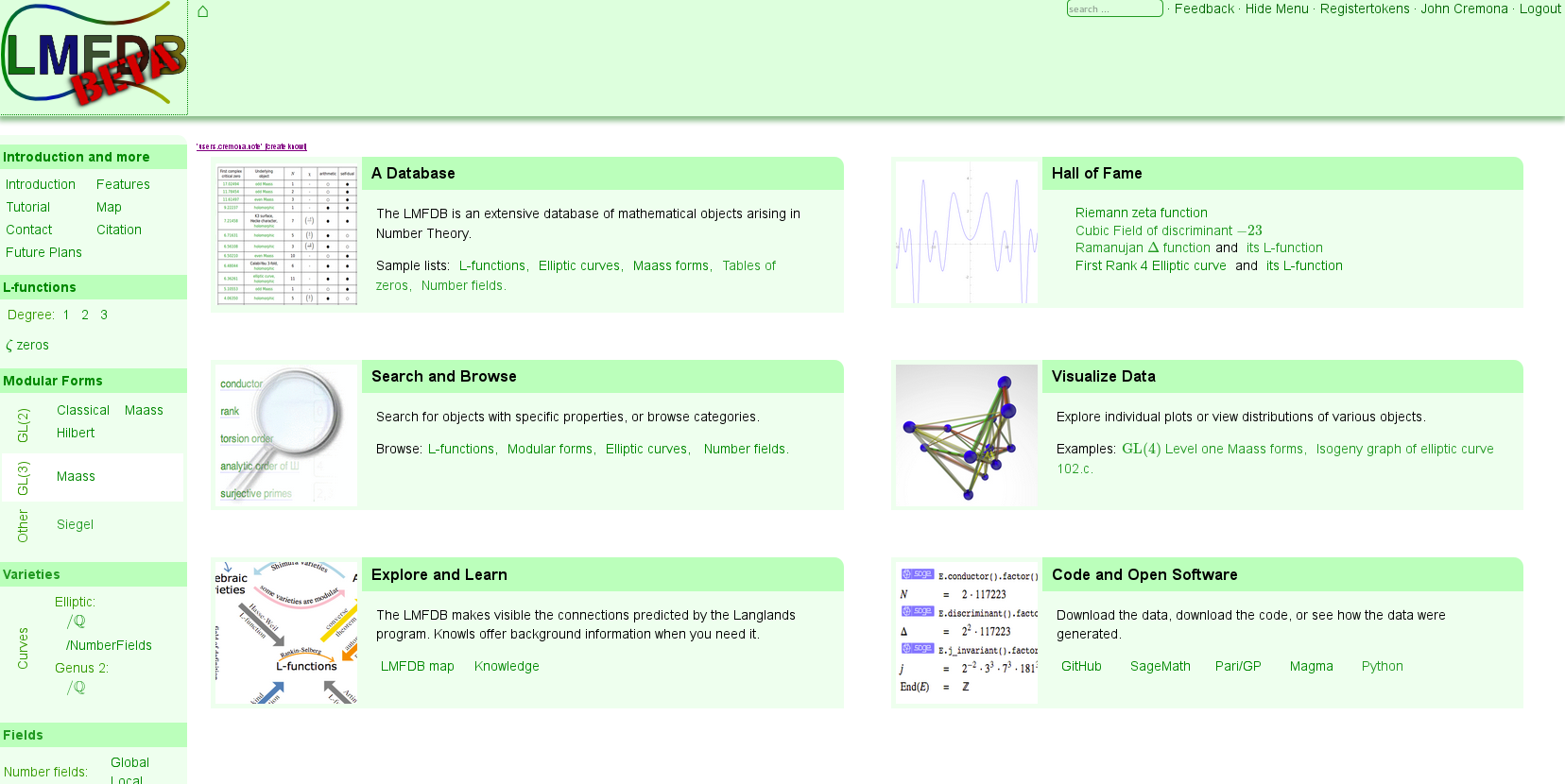}\hfill
}
\begin{center}The LMFDB home page at \url{www.lmfdb.org}, January 2016\end{center}

\section{L-functions and why they are important}
\textbf{L-functions} are at the heart of the LMFDB.  What are they?
We will give a brief survey, referring to number theory textbooks for
details.

The simplest L-function is the \textbf{Riemann zeta function}
$\zeta(s)$.  This
\begin{itemize}
\item is a \textit{complex analytic function} (apart from a pole at
  $s=1$); \item has a \textit{Dirichlet series} expansion over
  positive integers  (valid when  $\Re(s)>1$):
\[
   \zeta(s)=\sum_{n=1}^{\infty} \frac{1}{n^s};
\]
\item has an \textit{Euler product} expansion over primes~$p$ (when
  $\Re(s)>1$):
\[
  \zeta(s)=\prod_p(1-p^{-s})^{-1};
\]
\item satisfies a \textit{functional equation}:
\[
  \xi(s)=\pi^{-s/2}\Gamma(s/2)\zeta(s) = \xi(1-s);
\]
\item has links to the \textit{distribution of primes}.
\end{itemize}

\subsection{L-functions: a definition}
The definition of an L-function encapsulates these properties: it is a
\textbf{complex function} with a \textbf{Dirichlet series} and an
\textbf{Euler product} expansion which satisfies a \textbf{functional
  equation}.  There are other more technical axioms (by Selberg) which
we omit here: refer to the LMFDB's own knowledge database for details:
\url{http://www.lmfdb.org/knowledge/show/lfunction}.

Some of the defining properties have not in fact been proved for all the
types of L-function in the database: this can be very hard!  For example,
Andrew Wiles proved Fermat's Last Theorem by proving the modularity of
certain elliptic curves over $\Q$, which amounted to showing that the
L-functions associated to elliptic curves really are L-functions in
the above sense.  This is not yet known in general for elliptic curves
defined over other algebraic number fields.

Other expected properties of L-functions are not even known for
$\zeta(s)$.  For example, the \textbf{Riemann hypothesis} concerning the
the zeros of $\zeta(s)$
has remained open since it was formulated by Riemann in 1859.

\subsection{The Riemann Hypothesis}
The Riemann Hypothesis states that all the ``non-trivial'' zeros of
$\zeta(s)$ (excluding those coming trivially from poles of
$\Gamma(s)$) are on the ``critical line'' $\Re(s)=1/2$.

This was (part of) Hilbert's 8th problem and is also one of the Clay
Mathematics Institute Millennium Prize Problems, so a million dollars
awaits the person who proves it.  There are similar conjectures about
the location of the zeros of all L-functions, which are collectively
known as the \textbf{Generalized Riemann Hypothesis} (GRH).  These are
not only of theoretical (or financial!) interest, but have important
applications to the complexity of computing important quantities in
number theory.  For example,  computing the class number of
a number field is much faster if one assumes GRH for the number field's
own L-function, its \textbf{Dedekind $\zeta$-function}.

\medskip

What can a database say in relation to this problem?

\medskip

It can give the object its own web page
(\url{http://www.lmfdb.org/L/Riemann/}) which shows basic facts about
it, and its graph along the critical line $1/2+it$ to ``show'' the
first few zeroes.  This is a pedagogical function of the database.

It can also store all the zeroes which have so far been explicitly
computed: there are more than $10^{11}$ (that is \textit{one hundred
  billion}) of them at \url{http://www.lmfdb.org/zeros/zeta/}, all
computed to 100-bit precision by David Platt (Bristol), who in 2014
won a prize for his contributions to progress on the Goldbach
Conjecture.  This resource can then be used to study properties of the
zeroes, such as their distribution, and connections to random
matrices, showing that the database also serves as a research tool.

\subsection{Degrees of L-functions}
The \textbf{Euler product} for a general L-function has the form
\[
    L(s) = \prod_p 1/P_p(1/p^{s})
\]
where each $P_p(t)$ is a polynomial, and the product is over all
primes~$p$.  These polynomials all have the same degree, called the
\textbf{degree} of the L-function, except for a finite number indexed
by primes dividing an integer called the \textbf{conductor} of the
L-function, where the degree is smaller.  The zeros of these
polynomials are also restricted in a way depending on another
parameter, the \textbf{weight}.

For example, $\zeta(s)$ has $P_p(t)=1-t$ for all primes~$p$; the
degree is $d=1$, and the conductor is $N=1$.

\subsection{L-functions of degree $1$}
There are other L-functions of degree~$1$, with larger conductor~$N$,
which have been studied since the 19th century: \textbf{Dirichlet
  L-functions}.  Their Dirichlet coefficients $a_n$ are given by the
values of a \textbf{Dirichlet character} $a_n=\chi(n)$, meaning that
they are \textbf{multiplicative} and \textbf{periodic} with period
$N$.

An example with $N=4$ is
\[
  L(\chi,s)=1^{-s}-3^{-s}+5^{-s}-7^{-s}+-\dots,
\]
with all even coefficients~$0$ and the odd coefficients
alternating~$\pm1$.  Dirichlet used such L-functions to prove his
celebrated theorem about primes in arithmetic progressions: for any
integers $N\ge1$ and $a$, there are infinitely many primes $p\equiv
a\pmod{N}$, provided that $a$ and~$N$ are coprime. The previous
example can be used not only to show that there are infinitely many
primes $p\equiv1\pmod4$ (for which $\chi(p)=+1$) and infinitely many
primes $p\equiv3\pmod4$ (for which $\chi(p)=-1$) , but also to show
that (in a precise sense) the primes are \textit{equally distributed}
between these two classes.

This is a \textit{complete} list of all L-functions of degree~$1$.
For degrees greater than~$1$, a complete classification has not yet
been established, though a wide variety of sources of L-functions is
known, and in some cases (such as in degree~$2$, see below) we
conjecture that all L-functions do arise from these known sources.

\subsection{Other sources of  L-functions}
A wide variety of mathematical objects have L-functions: algebraic
number fields, algebraic varieties (including curves).  There is a
general term \textbf{motive} for objects which have L-functions.

In many cases, while we know how to define the L-function of a more
complicated object, it has not yet proved that it actually satisfies
the defining axioms for L-functions.  Even for elliptic curves
over~$\Q$, this would have been true until the mid-1990s; for elliptic
curves over real quadratic fields such as $\Q(\sqrt{2})$ it was true
until 2013!  Now, these elliptic curves are known to be modular
\cite{Siksek}.

\subsection{L-functions of number fields}
An \textbf{algebraic number field}, or simply \textbf{number field},
is a finite extension of the rational field $\Q$, such as
$\Q(\sqrt{2})$ or $\Q(i)$ or $\Q(e^{2\pi i/m})$.  Every number
field~$K$ has an L-function called its \textbf{Dedekind zeta function}
$\zeta_K(s)$, defined in a similar way to Riemann's
$\zeta(s)=\zeta_{\Q}(s)$, and with similar analytic properties.

Just as the analytic properties of $\zeta(s)$ imply facts about the
distribution of primes, from the analytic properties of $\zeta_K(s)$
we can deduce statements about prime factorizations in the field $K$.
For example, taking $K=\Q(e^{2\pi i/m})$ we can prove Dirichlet's
Theorem on primes in arithmetic progressions using a combination of
algebraic and analytic properties of~$\zeta_K(s)$.

Also, just as some properties of $\zeta(s)$ are not yet proved (for
example the Riemann Hypothesis), the same is true for $\zeta_K(s)$:
the Generalized Riemann Hypothesis or GRH remains unsolved.

\subsection{L-functions of curves}
Algebraic \textbf{curves} defined over algebraic number fields also
have L-functions, whose degree depends
on both the degree of the field over which the curve is defined, and
the genus of the curve. So an elliptic curve over $\Q$, which is a
curve of genus~$1$ defined over a field of degree~$1$, has a
degree~$2$ L-function, elliptic curves over fields of degree~$d$ have
L-functions of degree~$2d$, and so on.

It is widely believed that all degree~$2$ L-functions arise as
follows: they are either products of two degree~$1$ L-functions, or
come from elliptic curves over~$\Q$, or from (a special kind of)
\textbf{modular form}.  The insight of Weil, Taniyama, Shimura and
others in the 1960s and 1970s was to realize that the latter two
sources actually produce the \textit{same} L-functions!  This insight
is behind the famous theorem of Wiles \textit{et al.} that ``every
elliptic curve (over $\Q$) is modular'', from which Fermat's Last
Theorem was a consequence.  But it is still an unsolved problem to
show that those degree~$2$ L-functions which are not products of
degree~$1$ L-functions do all arise from automorphic forms.

\subsection{Higher degree L-functions}
For degrees~$3$ and~$4$, we do not yet even have a conjecture
concerning all sources of L-functions, and for those which are known,
not all the conjectured connections between them have been proved.

We mentioned above the recent result \cite{Siksek} that elliptic
curves defined over real quadratic fields (such as $\Q(\sqrt{5})$) are
modular.  This means that two sources of L-functions of degree~$4$:
on the one hand, elliptic curves over such a field, and on the other
hand \textbf{Hilbert modular forms} over the same field, actually
produce the \textit{same} L-functions.  Such results are extremely
deep and require a vast amount of theory to establish, including real,
complex and $p$-adic analysis and algebra, as well as some explicit
computations (the ArXiV version of \cite{Siksek} includes a number of
\Magma\ scripts).

By contrast, over \textit{imaginary} quadratic fields
(e.g. $\Q(\sqrt{-1})$) we conjecture, but cannot prove in general,
that elliptic curves have L-functions also attached to a different
kind of modular form, \textbf{Bianchi modular forms}.  These can be
computed, and work is in progress in entering many examples into the
LMFDB, even though they are not all known to ``be modular'' and hence
have genuine L-functions.

Modularity of \textit{individual} elliptic curves over imaginary
quadratic fields can be proved using the Serre-Faltings-Livn\'e method
(which uses Galois representations rather than analysis) as explained
in a 2008 paper \cite{DGP} by Dieulefait, Guerberoff and Pacetti.  We
are currently using their method to prove modularity of \textit{all}
the curves in the database; at the same time we are developing
enhancements to the algorithm to make it more efficient.  A
theoretical proof that all elliptic curves over these fields are
modular seems very far off, so even in the world of L-functions of
degree~$4$ it is still important to carry out experiments and collect
data.

\subsection{Showing connections through the LMFDB}
The LMFDB shows connections between different objects with the same
L-function, such as those described above, by linking its databases of
(for example) elliptic curves over real quadratic fields, and Hilbert
modular forms over the same field.  The home page of each elliptic
curve includes a link to the associated Hilbert modular form, and to
the associated L-function, and (in progress) vice versa.

One difficulty we have encountered in setting up these links on the
website, which is perhaps typical in a large project where many
different individuals are providing data, is to maintain
\textit{consistency of labelling} of objects.  Over the field
$\Q(\sqrt{5})$, the Hilbert modular forms were computed (in \Magma) by
John Voight (Dartmouth College) and Steve Donnelly (Sydney) \cite{DV},
while the elliptic curves were computed (in \Sage) by Jonathan Bober
(Bristol), William Stein (Washington), Alyson Deines (CCR) and others
\cite{Bober}.  These groups used essentially the same naming
convention, but we were careful to check that the labels of matching
objects did match exactly, resulting in one set of data (the elliptic
curves) requiring relabelling.

\section{The LMFDB database}
The LMFDB consists of both a database, where the data collection
itself is organized and stored, together with the website
\url{www.lmfdb.org}.  This provides a sophisticated user interface to
the data, has home pages for individual objects in the database,
showing links between related objects, and also provides an online
repository of knowledge about L-functions and related objects, through
its \textbf{knowledge database}.

Both database and website are currently hosted on servers at Warwick,
funded by EPSRC; until 2013 they were hosted at the University of
Washington on NSF-funded servers administered by William Stein.  Plans
are also underway to have mirror sites in other countries: this is an
international project.

The LMFDB is also a group of mathematicians who collaborate to create
and develop the database and its website.  We will say more about this
collaboration in the final section.

\subsection{The database and website software}
We are using the open-source database software \mongodb.  This
currently holds nearly a terabyte of data and indices.  This choice
was made because \mongodb\ allows data to be organized in a completely
flexible schema, rather than having to specify the schema for each
item in advance as with SQL databases.  It also has a powerful
\Python\ interface, \pymongo, which suits the project well, since it
allows the website code to use other \Python\ modules such as
\flask\ (a web framework), and to have access to all the power of
\Sage, another large \Python-based open source mathematical software
project.  All of these are open source, which is another essential
requirement.  (Note, however, that not all of the data in the database
have been computed using open source tools.)  The website code is a
collaborative open source project hosted at
\GitHub\ (see \url{https://github.com/LMFDB/lmfdb}).

Basing all the website code on \Python\ has many advantages.  It is
relatively easy to learn to use, which is important since we want the
barriers to new people joining our project to be as low as possible.
And it is phenomenally powerful, giving access to a vast array of
additional modules for interfacing with the database (\pymongo),
running the web framework (\flask), web page templating (\jinja),
testing and more.

Anyone contributing to the project who wants to do more than just
donate data has to learn how to use this software.  At project
workshops we run tutorial sessions for newcomers, where code is
written by beginners under the guidance of more experienced peers.
All code is reviewed and tested before being adopted, as well as being
subject to some automated testing.

\subsection{Database organization}
The database as a whole consists of around 35 individual databases
containing collections of mathematical objects (including {\tt
  elliptic\_curves}, {\tt hilbert\_modular\_forms}, and {\tt
  number\_fields}) and other data such as the {\tt knowledge}
database, which holds the contents of \textbf{knowls} (see below).

The data are indexed in various ways for faster searching, and, of
course, backed up regularly.  Many parts of the database can also be
recreated from plain text data files which are stored in separate
\git\ repositories, also hosted on \GitHub.

Each constituent database contains collections of records,
and these records hold the data in a flexible format: additional data
fields can be added later.

\subsection{Sample database entry}
To take just one example, the database {\tt number\_fields} contains
just one collection {\tt fields}, for which a typical entry looks as
follows (after being converted by \pymongo\ into a \Python\
dictionary):
\begin{verbatim}
{u'_id': ObjectId('4cb80fdb5009fb52db0946b6'),
 u'class_group': u'',
 u'class_number': 1,
 u'coeffs': u'1,0,-1,1',
 u'degree': 3,
 u'disc_abs_key': u'00123',
 u'disc_sign': -1,
 u'galois': {u'n': 3, u't': 2},
 u'label': u'3.1.23.1',
 u'ramps': [u'23'],
 u'signature': u'1,1',
 u'unitsGmodule': [[3, 1]]}
\end{verbatim}
Here we see the coefficients of the minimal polynomial $x^3-x^2+1$ of
a generator of the field stored as {\tt coeffs}, and the label
'3.1.23.1' which also uniquely determines the field.  Invariants of
the field which are easy to compute on the fly, and to which we do not
need to provide direct access through database queries, need not be
stored, while quantities which might be expensive to compute, or for
which we may want to run searches, are stored and indexed.

This is only a simple example.  The database entry for an individual
elliptic curve over~$\Q$ currently contains 33 fields, some very
technical.  The number of fields grows over time as new data are
contributed.  For example, in 2015 Jeremy Rouse offered to provide
information concerning the $2$-adic Galois Representation attached to
every elliptic curve over~$\Q$, after developing and implementing an
algorithm to determine this jointly with David Zureick-Brown (see
\cite{Rouse}).  He provided us with a \Magma\ script of their
implementation, we ran it and uploaded the data, and added a
corresponding section on the home page of every curve showing these
additional data.

\subsection{Software choices, pros and cons}
Using off-the-shelf software has plenty of advantages but will never
be perfect for a mathematical project.

Most mathematicians, even those with substantial computational
experience and expertise, know almost nothing about databases or
running websites, and many of the contributors to the LMFDB knew
nothing at all about these before they joined the project.  Decisions
about the specific software used by the project was made by those who
did have such experience, notably William Stein (lead developer of
\Sage) and Harald Schilly (another key developer of the {\sc
  SageMathCloud} project, \url{https://cloud.sagemath.com/}).

We have already seen some of the advantages of our choice of database,
\mongodb.  There are disadvantages too: \mongodb\ data consists of
strings or integers or floating point values, with strings as keys.
Values can also be lists of these, but a serious deficiency for
number-theoretic data is that the integers cannot be larger
than~$2^{32}$.  This means that most data fields which hold integers
have to be stored as strings, and this limits functionality, such as
searching for the value being in a certain range.

Similarly, rational numbers cannot be stored as such, or even as a
pair of integers [\textit{numerator},\textit{denominator}] if these
could be large, so instead they are stored as strings such as `{\tt
  1728}' or `{\tt -122023936/161051}'.  Building on these,
considerable thought has to be given as to how to store more
complicated data, such as an element of a number field.  Decisions
such as these are made by consensus at LMFDB workshops, since they
affect all developers, even though the effects of such decisions are
hidden from users of the website.

\section{The LMFDB website}
The LMFDB website serves several purposes.  It provides
\begin{itemize}
\item a shop window for the data;
\item a way to visualize the data, and the connections between
  different, linked mathematical objects;
\item a way to browse types of object;
\item a way to search for objects with specified properties;
\item a repository of knowledge through its ``knowledge database'';
\item a source of data for downloading for further work.
\end{itemize}

Catering for several different audiences at once is hard to
get right!

\subsection{Technical support (or lack of)}
The project would benefit greatly from having technical support staff.
Our current grant from EPSRC does not provide this---it does support
six postdoctoral researchers, who all have a certain amount of
experience writing mathematical software, but not any dedicated
software engineers.  For this, we are currently relying on charitable
contributions of time.  We would not be where we are now, and indeed
the website would never have been launched, without the enormous
contributions of one person in particular: Harald Schilly%
, a doctoral
student in Vienna and software consultant, who knows more than the
rest put together about \Python, \mongodb, \flask, and the rest.

From September 2015, through the Horizon 2020 European Research
Infrastructure project OpenDreamKit (\url{http://opendreamkit.org/}),
which provides substantial funding for the development of open source
computational mathematics, we are currently seeking to employ a
software engineer to provide support to the project.

\subsection{Homepages}
A key organizing principle of the LMFDB is that every object has its
own homepage.  These have mathematically meaningful, permanent URLs
which follow a carefully thought out schema.  The home pages
themselves are created on demand from templates, filled in with data
partly retrieved directly from the database and partly computed on the
fly.  For example, the elliptic curve with label {\tt 5077a1} has URL
\url{http://www.lmfdb.org/EllipticCurve/Q/5077/a/1}.

Each homepage gives a view of the object (depending on its nature),
highlighting its most important properties, with \textit{breadcrumbs}
to show its position in the whole.  In some cases, where the object
has some interesting additional historical or mathematical
significance, this can also be shown on its home page.  For example,
the elliptic curve {\tt 5077a1} was used by Dorian Goldfeld in 1985 to
solve Gauss's class number problem effectively, by making use of a new
connection between the problem and L-functions of elliptic curves, and
this piece of historical information is shown on the curve's home
page.

A \textbf{related objects} box on each homepage provides links between
related objects.  For example, from the pages of an elliptic curve, or
a number field, or a modular form there are links to the associated
L-function.

Where possible, on the home page of an object, we make it possible to
see and download code which will re-create the object in one of the
standard number-theoretical packages (\Sage, \PARI\ or \Magma) and
work with it there.  In this way, the LMFDB can be used by students
learning a subject who wish to work out their own examples, as well as
researchers wishing to carry out larger-scale investigation starting
from the LMFDB data.  A more sophisticated programming interface
through \Sage\ is also planned.

\subsection{Searching and browsing}
Each class of objects in the LMFDB has its own \textbf{Browse and
  Search} page.

The \textbf{Browse} section is intended to be usable by people who know nothing
of the underlying theory but want to browse through examples without
having to type anything or have technical knowledge.

The \textbf{Search} section is more for experts looking for a specific
object (possibly by its label), or for an object with certain
properties: {``a number field with Galois group $C_5$ ramified
  only at $p=5$''}, or {``an elliptic curve with rank $2$ and
  non-trivial Tate-Shafarevich group''}, or {``a classical
  modular form of weight~$12$ and level $12$''}.  This leads to a
\textbf{Search Results} page listing all database entries which match
(if any), with links to the home pages of each individual matching
object.

\subsection{Knowledge and knowls}
The knowledge aspect of the LMFDB exists in the first place as a
glossary of technical terms used on the web pages, so the pages
themselves do not get cluttered up, and there is consistency between
pages on basic definitions.

The mechanism which serves these is the \textbf{knowl}, created by
Harald Schilly and first demonstrated at an LMFDB workshop.  The text
expands within the page and can be dismissed after reading, without
any need for ``pop-ups'' or new pages.

Knowls can be used anywhere on the web---for example, I use them on my
own web page of preprints and publications to display abstracts of
papers.
Another good example of their use is in the online undergraduate
textbook on Linear Algebra by Robert Beezer \cite{Beezer}.  For more
about knowls and how to use them, see the knowl on the LMFDB itself,
\url{http://www.lmfdb.org/knowledge/show/doc.knowl}, or the page
\url{http://aimath.org/knowlepedia/}.

The content of knowls can be edited by any project member (someone who
has a login account), and is itself stored in the database.

\section{The LMFDB project}

\subsection{The LMFDB as a collaborative project}
The LMFDB was first conceived at an AIM workshop in 2007. It holds
regular workshops, which are run
along the lines of AIM workshops: few talks and a lot of hard work.
As well as individual workshops of around 30 people, there are smaller
groups who meet to work together on specific projects, and there have
also been longer periods of activity hosted at MSRI, during the
semester programme ``Arithmetic Statistics'' in 2011, and ICERM,
during the semester programme ``Computational Aspects of the Langlands
Program'' in late 2015 (see \url{http://icerm.brown.edu/sp-f15/}).
All members of the organizing committee for the latter are LMFDB
contributors, and we expect that the LMFDB will make a substantial
leap forward during the semester.

The AIM connection remains strong: both Brian Conrey (Director of AIM)
and David Farmer (Director of Programs at AIM) are number theorists
who have been intimately connected with the project from the start.

We have Editorial and Management Boards, but essentially all
decisions are made by consensus at workshops.

\subsection{Funding}
During 2008--2012, the LMFDB was funded by NSF FRG grant DMS:0757627;
currently (2013-2019) it is supported by Programme Grant EP/K034383/1
from the UK research council EPSRC.  The investigators on this are the
author and Samir Siksek (Warwick) and Brian Conrey (AIM and Bristol),
and Andy Booker and Jon Keating (Bristol).  David Farmer (AIM) is a
project partner, as are Fernando Rodriguez-Villegas (ICTP), William
Stein (Washington) and Mike Rubinstein (Waterloo).

These research grants provide funding both for LMFDB workshops and for
servers hosting the database and website; the NSF FRG grant also paid for
some technical software support.

\subsection{Collaboration}
The LMFDB encompasses such a wide range of mathematics, and it is
essential to have an equally wide range of mathematical expertise
contributing to the project.  Many of the collaborators on the LMFDB
project, who are all listed at
\url{http://www.lmfdb.org/acknowledgment}, have contributed not by
coding for the website but by providing the data (without which the
project would be nothing!).  More contributors are always welcome.

\section{Acknowledgements}
Funding for the LMFDB project as a whole has been mentioned above.
For a full list of contributors to the LMFDB project, financial and
other support, see \url{http://www.lmfdb.org/acknowledgment}.  The
author personally acknowledges support from EPSRC Programme Grant
EP/K034383/1 ``LMF: L-Functions and Modular Forms'' and (since 1
September 2015) from the OpenDreamKit Horizon 2020 European Research
Infrastructures project (\#676541).

\end{document}